\documentclass[12pt]{article}
\usepackage{fullpage,amsmath,amssymb,graphicx}

\font\cyr=wncyr10

\newtheorem{thm}{Theorem}[section]

\newtheorem{rem}[thm]{Remark}

\newcommand{\Ga}{\Gamma}

\newcommand{\Sa}{\Sigma}

\newcommand{\Da}{\Delta}
\newcommand{\e}{\epsilon}

\newcommand{\la}{\lambda}
\newcommand{\La}{\Lambda}

\newcommand{\Om}{\Omega}

\newcommand{\vp}{\varphi}

\newcommand{\sch}{{\cal H}}
\newcommand{\sci}{{\cal I}}

\newcommand{\bc}{{\mathbb C}}
\newcommand{\C}{{\mathbb C}}
\newcommand{\bF}{{\mathbb F}}

\newcommand{\bh}{{\mathbb H}}

\newcommand{\bo}{{\mathbb O}}
\newcommand{\bp}{{\mathbb P}}

\newcommand{\br}{{\mathbb R}}
\newcommand{\R}{{\mathbb R}}

\newcommand{\gtH}{{\mathfrak H}}

\newcommand{\gtS}{{\mathfrak S}}

\newcommand{\p}{\partial}

\newcommand{\ra}{\rightarrow}
\newcommand{\lra}{\longrightarrow}

\newcommand{\bs}{\backslash}
\newcommand{\ov}{\overline}
\newcommand{\col}{\!:\!}

\newcommand{\fh}{H_{\bF}^}
\newcommand{\hh}{H_{\bh}^}
\newcommand{\rh}{H_{\br}^}
\newcommand{\oh}{H_{\bo}^2}

\newcommand{\vol}{\operatorname{vol}}
\newcommand{\orb}{\operatorname{orb}}

\newcommand{\rint}{\text{int}}

\newcommand{\Sp}{\operatorname{Sp}}

\newcommand{\SO}{\operatorname{SO}}

\newcommand{\F}{\operatorname{F}}
\newcommand{\SU}{\operatorname{SU}}
\newcommand{\Spin}{\operatorname{Spin}}
\newcommand{\id}{\operatorname{id}}

\newcommand{\im}{\operatorname{Im}}

\newcommand{\isom}{\operatorname{Isom}}

\newenvironment{pf}{\begin{trivlist}\item[]{\bf Proof:\ }}
{\mbox{}\hfill\rule{.08in}{.08in}\end{trivlist}}

\begin{document}
\title{Rigidity of locally symmetric rank one manifolds of infinite volume}
\author{Boris N. Apanasov}
\date{}
\maketitle

\begin{abstract} We discuss questions by
Mostow \cite{Mo1}, Bers \cite{B} and Krushkal \cite{Kr1, Kr2} about uniqueness of a conformal or spherical CR structure on the sphere at infinity $\partial H_\mathbb{F}^n$ of symmetric rank one space $H_\mathbb{F}^n$ over division algebra $\bF=\br\,,\bc\,,\bh\,,\text{or}\,\, \bo $ compatible with the action of a discrete group
$G\subset\operatorname{Isom}H_\mathbb{F}^n$. Introducing a nilpotent Sierpi\'{n}ski carpet with a positive Lebesgue measure in the nilpotent geometry in
$\partial H_\mathbb{F}^n\setminus\{\infty\}$ and its stretching, we construct a non-rigid discrete $\mathbb{F}$-hyperbolic groups $G\subset\operatorname{Isom}H_\mathbb{F}^n$ whose non-trivial deformations are induced by $G$-equivariant homeomorphisms. Here we consider two situations: either the limit set $\Lambda(G)$ is the whole sphere at infinity $\partial H_\mathbb{F}^n$ or restrictions of such non-trivial deformations to components of the discontinuity set $\Om(G)\subset \partial H_\mathbb{F}^n$ are given by restrictions of $\bF$-hyperbolic isometries. In both cases the demonstrated non-rigidity is related to non-ergodic dynamics of the discrete group action on the limit set which could be the whole sphere at infinity.
\end{abstract}
    \footnotetext[2]{
\hfil\hfil\hfil}

    \section{Introduction}
Many authors study deformations of negatively curved locally symmetric spaces of $\R$-rank one, that is spaces modeled
on $\bF$-hyperbolic spaces $H^n_{\bF}$ over division algebras $\bF$ that are either
real $\br$,  complex $\bc$, quaternions $\bh$, or
Cayley numbers (octonions) $\bo$. Here the central result is the remarkable strong rigidity (see Mostow \cite{Mo1, Mo2}) claiming that homotopy equivalence of such finite volume spaces of real dimension at least 3 implies their isometry. It can be reformulated for discrete isometric action of their fundamental groups in $\bF$-hyperbolic spaces:

\begin{thm}\label{Mostow} Let $\Ga_1, \Ga_2\subset\isom\fh m$ be isomorphic lattices in
$\bF$-hyperbolic space $\fh m$ with real dimension $n\geq 3$, that is isomorphic discrete isometry
groups with finite volume quotients $\fh m/\Ga_i<\infty$ over real $\R$, complex $\C$, quaternionic $\bh$ or Cayley $\bo$ numbers.  Then they are conjugate in
$\isom\fh m$, and the original isomorphism is induced by this conjugation.
\end{thm}

Examining the main stages of the proof of this rigidity theorem shows that the main question in extending this rigidity to a wider class of
$\F$-hyperbolic structures which may have infinite volume, is concerned with
the uniqueness of a spherical CR-structure (or conformal one for real $\F=\R$) on
the sphere at infinity $\p\fh m$
compatible with the action of a discrete group $\Ga\subset\isom\fh m$.  It was
questioned by Mostow \cite{Mo1} and Krushkal \cite{Kr1, Kr2} whether such uniqueness
takes place in the case of those discrete $\F$-hyperbolic groups whose limit sets
$\La(\Ga)$ coincide with the whole sphere $\p\fh m$.

The goal of this paper is to address these deformation questions and to show that in general, the uniqueness of such a spherical CR-structure (or conformal one for real $\F=\R$, see \cite{A1, A3}) at infinity $\p\fh m$
and the rigidity of the corresponding $\bF$-hyperbolic structure fails:

\begin{thm}\label{non-rigid}
Let $\fh n$ be a symmetric $\R$-rank one space over division algebra $\bF$ (either $\R$, $\C$, $\bh$, or $\bo$) of real dimension at least 3 whose nilpotent Carnot group $\gtH\subset\isom\fh n$ acting at infinity
$\p\fh n\bs\{\infty\}$ has non-rigid lattices. Then there are non-rigid discrete $\bF$-hyperbolic groups $\Ga\subset\isom\fh n$ whose limit set $\La(\Ga)$ is the whole sphere $\p\fh n$.
\end{thm}

Especially interesting case is related to quaternionic hyperbolic spaces $\hh n$ or the Cayley (octonionic) hyperbolic plane $\oh$. These geometries are symmetric rank one spaces having nevertheless the super rigidity of deformations of their lattices, see Corlette \cite{C}, Pansu \cite{P} , Gromov-Schoen \cite{GS} and Kamishima \cite{Ky}.

Our construction of non-rigid discrete $\bF$-hyperbolic groups $\Ga\subset\isom\fh n$ in the proof of this result allows some extension to non-trivial deformations of negatively curved locally symmetric $\R$-rank one spaces with boundary components at infinity where restrictions of homeomorphisms inducing such deformations are extensions of $\bF$-hyperbolic isometries to those boundary components. This answers another question
of L.~Bers \cite{B} and S.~Krushkal \cite{Kr1} and its analogue in $\bF$-hyperbolic spaces on the support of deforming homeomorphisms at the sphere at infinity $\p\fh n$:

\begin{thm}\label{CR-boundary} For discrete $\bF$-hyperbolic groups $\Ga\subset\isom\fh n$ isometrically acting in
symmetric $\R$-rank one spaces in Theorem \ref{non-rigid} with non-empty discontinuity set
$\Om(\Ga)\subset\p\fh n$,
there are their non-trivial deformations induced by continuous family of $\Ga$-equivariant homeomorphisms
$\{\Phi^0_t\}$  whose restrictions to connected components of $\Om(\Ga)$ are restrictions of some $\bF$-hyperbolic isometries.
In particular, for real hyperbolic spaces $\rh n$, $n=2,3$, there are Beltrami $\Ga$-differentials $\mu(z)dz/d\ov{z}$ whose supports are contained in the limit set $\La(\Ga)\subset\ov\C$ (of positive Lebesgue (n-1)-measure).
\end{thm}

We should note that in both Theorems \ref{non-rigid} and \ref{CR-boundary}, the crucial property of our constructed non-rigid discrete $\bF$-hyperbolic groups $\Ga$ is that the limit
subset $\La(\Ga)\cap\ov{P(\Ga)}$ on the boundary of the fundamental polyhedron
$P(\Ga)\subset\fh m$  has a positive Lebesgue measure in the sphere at infinity $\p\fh n$. This is closely related
to non-ergodic dynamics of our $\bF$-hyperbolic discrete holonomy group $\Ga$ action on the limit set  $\La(\Ga)$ in the sphere at infinity  $\p\fh n$ of the $\bF$-hyperbolic space $\fh n$, see Remarks \ref{dynamics}, \ref{measure}, \ref{convex hull} and \cite{S, A3, A6}).

\section{Preliminaries}

\subsection{Symmetric spaces of rank one}

The symmetric spaces of $\br$-rank one of non-compact type
are the
hyperbolic spaces $\fh n$ over four division algebras $\bF$ - either
the real numbers $\br$, or
the complex numbers $\bc$, or
the quaternions $\bh$, or the Cayley numbers (octonionions) $\bo$; in last case $n=2$.
They are respectively called as real, complex, quaternionic and octonionic hyperbolic
spaces (the latter one $\oh$ is also known as the Cayley hyperbolic plane).
Algebraically these spaces can be described as the corresponding quotients:
$\SO(n,1)/\SO(n)$,
$\SU(n,1)/\SU(n)$,
$\Sp(n,1)/\Sp(n)\cdot\Sp(1)$ and
$\F_4^{-20}/\Spin (9)$
where the latter group $\F_4^{-20}$ of automorphisms of the Cayley plane $\oh$
is the real form of $\F_4$ of rank one. We normalize
the metric so the (negative) sectional curvature of $\fh n$ is bounded
from below by $-1$.

Following Mostow \cite{Mo1} and using the standard involution (conjugation) in $\bF$,
$z\to\bar z$, one can define projective models of the hyperbolic spaces
$\fh n$ as the set of negative lines in the Hermitian vector space $\bF ^{n,1}$,
with Hermitian structure given by the indefinite $(n,1)$-form
$$
\langle\langle z,w\rangle\rangle
=z_1\overline w_1+\cdots+z_n\overline w_n-z_{n+1}\overline w_{n+1}\,.
$$
The boundary $\p\fh n$ consists of all null $\bF$-lines,
$$
\p \fh n = \{[z]\in \bp\bF ^{n,1}\col \langle\langle z,z\rangle\rangle=0\}\,,
$$
 and is homeomorphic to the sphere $S^{kn-1}$ with $k=\dim_{\R}\bF$.

Taking non-homogeneous coordinates, one can obtain ball models (in the unit ball
$B^n_{\bF}$ in $\bF^n$) for the first three spaces. Here we note that since the multiplication by quaternions
is not commutative, we specify that we use
``left" vector space ${\mathbb H}^{n,1}$
where the multiplication by quaternion numbers is on the left.
However, it does not work for
the Cayley hyperbolic plane since $\bo$ is non-associative, and one should use a Jordan algebra of
$3\times 3$ Hermitian matrices with entries from $\bo$ whose group of automorphisms is
$\F_4$, see \cite{SV, Mo1}.

\subsection{Upper half-space model of $\fh n$}

Another models of $\fh n$ use the so called horospherical coordinates \cite{AX, AK, G}
based on foliations
of  $\fh n$ by horospheres centered at a fixed point $\infty$ at infinity $\p\fh n$ which is
homeomorphic to $(n\dim_{\br}\bF-1)$-dimensional sphere. Such a horosphere can be
identified with the nilpotent
group $N$ in the Iwasawa decomposition $KAN$ of the automorphism group
of $\fh n$. The nilpotent
group $N$ can be identified with the product $\gtH=\bF^{n-1}\times\im \bF$ (see \cite{Mo1})
equipped with the operations:
$$
(\xi,v)\cdot (\xi',v')=(\xi+\xi',v+v'+2\im \langle \xi,\xi'
\rangle )\quad \text{and}\quad (\xi,v)^{-1}=(-\xi,-v)\,,
$$
where $\langle ,\rangle$ is the standard Hermitian product
in $\bF^{n-1}$, $\langle z,w\rangle =\sum z_i\ov{w_i}$.
The group $\gtH$ is a 2-step
nilpotent Carnot group with center $\{0\}\times \im\bF \subset \gtH=\bF ^{n-1}\times \im\bF$,
and acts on itself by the left translations $T_h(g)=h\cdot g\,,\quad
h,g\in \gtH$.

Now we may identify
$$
\fh n\cup \p \fh n\bs\{\infty\} \lra \gtH \times [0,\infty)=
\bF^{n-1}\times \im\bF\times [0,\infty)\,,
$$
and call this identification the {\it ``upper half-space model"} for $\fh n$
with the natural horospherical coordinates $(\xi,v,u)$.
In these coordinates, the above left action
of $\gtH$ on itself extends to
an isometric action (Carnot translations) on the $\bF$-hyperbolic space
in the following form:
$$T_{(\xi_0,v_0)} \col (\xi,v,u)\longmapsto (\xi_0+\xi\,,v_0+v+2\im
\langle \xi_0,\xi\rangle \,,u)\,,
$$
where $(\xi,v,u)\in \bF^{n-1}\times \im\bF\times [0,\infty)$.

There are a natural norm and an induced by this norm distance on the Carnot group
$\gtH=\bF^{n-1}\times \im\bF$, which are known in the case of the
Heisenberg group (when $\bF=\bc$) as the Kor\'{a}nyi-Cygan's norm and distance.
  Using horospherical coordinates, they can be extended to
a norm on $\fh n$, see \cite {AX, AK}:
\begin{eqnarray}\label{norm}
|(\xi,v,u)|_c=|\,(|\xi|^2 + u - v)|_{\bF}^{1/2}\,,
\end{eqnarray}
where $|.|_{\bF}$ is the norm in $\bF$, and to
a metric $\rho_c$ on
the upper half-space model $\bF^{n-1}\times \im\bF\times (0,\infty)$ of
$H^n_{\bF}$:

\begin{eqnarray}\label{metric}
\rho_c\bigl((\xi,v,u),(\xi',v',u')\bigr)
=\bigl|\,|\xi-\xi'|^2 + |u-u'| -(v - v'+
2\im \langle\xi,\xi'\rangle )\bigr|_{\bF}^{\frac{1}{2}}\,.
\end{eqnarray}

It follows directly from the definition that
Carnot translations and  rotations are isometries
with respect to the Kor\'{a}nyi-Cygan metric $\rho_c$.
Moreover, the restrictions of this metric to different
horospheres centered at $\infty$ are the same, so the
Kor\'{a}nyi-Cygan metric plays the same role as Euclidean metric does on the upper half-space model
for the real hyperbolic space $\rh n$.

Non-zero numbers $\la\in \bF^*$ act on the Carnot group  $\gtH=\bF^{n-1}\times\im\bF\times\{0\}$
and on the complex hyperbolic space $\fh n=\bF^{n-1}\times\im\bF\times\R_+$ by Carnot dilations (which, for
$|\la|_{\bF}=1$, are in fact rotations):
\begin{eqnarray}\label{dilation}
D_{\lambda}\col (\xi,v,u) \longmapsto (\lambda \xi\,,|\lambda|_{\bF}^2v\,,|\lambda|_{\bF}^2u)\,.
\end{eqnarray}

\subsection{Spinal spheres in $\p\fh n$, bisectors and inversions}

For a $\bF$-hyperbolic hypersubspace $V=\fh {n-1}\subset\fh n$, its boundary at infinity $\p V\subset \p\fh n$ is a $\bF$-hyperchain. Assuming $\fh n=\F^{n-1}\times \im\F\times \R_+$, one has infinite hyperchains containing $\infty$ and finite ones. Let us denote $\bF$-inversion of order 2 with respect to $V$ by $\sci_V$. Its set of fixed points is $V$. If a hyperchain $\p V\subset\bF^{n-1}\times \im\bF\times \{0\}\cong\gtH$ is finite, its center is $\sci_V(\infty)\in\gtH$.
In the unit ball model $B_{\bF}^n\subset\bF^n$ of the $\bF$-hyperbolic space $\fh n$, the $\bF$-inversion $\sci_V$ acts as $\sci\,:\, (z',z_n)\mapsto (z',-z_n)$. Considering an isometry of this ball model to the upper half-space model of the $\bF$-hyperbolic space $\fh n$ (where one has
the nilpotent geometry $\gtH=\bF^{n-1}\times\im\bF\times\{0\}$ at infinity of $\ov{\fh n}\bs\{\infty\}$),
one has the hyperchain $\p V=(z',0)\in\p B^n_{\bF}$ as the hyperchain
$(z,\{0\})\in\bF^{n-1}\times\im\bF\col\,\, \mid z\mid=1$ in the upper half-space model of the $\F$-hyperbolic space
$\fh n$. Then the above
inversion $\sci$ in the unit ball is conjugate to the inversion in the finite hyperchain $\p V=\p\fh {n-1}$ in the upper half-space model of the $\F$-hyperbolic space $\fh n$ which acts in the Carnot group $\gtH$ as follows:
\begin{eqnarray}\label{inversion}
\sci(\xi,v)=\left (\frac {\xi}{|\xi|^2-v}\,,\, \frac {-v}{|v|^2+|\xi|^4}
\right )\quad \text {where}\,\,(\xi,v)\in \gtH=\F^{n-1}\times \im\F \,.
\end{eqnarray}

Obviously, the inversion $\sci$ interchanges the origin in the Carnot group $\gtH$ and
$\infty$ and preserves the unit sphere (in Kor\'{a}nyi-Cygan metric $\rho_c$):
\begin{eqnarray}\label{sphere}
S(0,1)=\{(\xi,v)\col\,\, |v|^2+|\xi|^4=1\}
\end{eqnarray}
containig the (bounded, horizontal)  hyperchain $\p V$ which is pointwise fixed
by $\sci$, and interchanges the points $(0,v)$ and $(0,-v)$, where $(0,v)$ are points in the unit sphere $S(0,1)$ which belong to the center $\{0\}\times\im\bF$ of the Carnot group $\gtH=\bF^{n-1}\times \im\bF$. In the case $\bF=\C$,
the Carnot group $\gtH$ is the Heisenberg group $\sch_{2n-1}=\C^{n-1}\times\R$, and
the inversion $\sci$ interchanges the unit sphere $S(0,1)$ poles $(0,1)$ and $(0,-1)$ in the center $\{0\}\times\R$ of the Heisenberg group $\sch_{2n-1}$.

For a given Carnot ball $B_i\subset \gtH$ at infinity of $\bF$-hyperbolic space $\fh n$, one has an isometry
$h_i\in\isom\fh n$ that maps this ball $B_i$ to the unit $\gtH$-ball $B(0,1)\subset\gtH=\bF^{n-1}\times\im\bF$
centered at the origin.
Then the $\bF$-hyperbolic isometry
\begin{eqnarray}\label{F-inversion}
\sci_i=h_i^{-1}\sci h_i\in\isom\fh n
\end{eqnarray}
is the inversion of order 2 that preserves
the boundary $\gtH$-sphere $S_i=\p B_i$. It interchanges the balls in $\ov{\gtH}$
bounded by this sphere and pointwise fixes the (horizontal) hyperchain
$c_i\subset S_i$.

The above $\gtH$-spheres $S_i=\p B_i$ in the Carnot group $\gtH=\bF^{n-1}\times\im\bF$ and their isometric $\bF$-hyperbolic images in the sphere at infinity $\p\fh n$ (sometimes called spinal spheres, see \cite{Mo1, G}) are boundaries of $\bF$-hyperbolic bisectors with respect the $\bF$-hyperbolic metric $d$ and pairs of points $z_1$, $z_2$ :
\begin{eqnarray}\label{bisector}
\gtS (z_1, z_2)=\{z\in \fh n
\col\,\, d(z_1, z)=d(z_2, z)\}\,
\end{eqnarray}
which are hypersurfaces in $\fh n$. In the case of $\bF=\R$, i.e. in the real hyperbolic space $\rh n$, they are totally geodesic hyperbolic hyperplanes. Bisectors in $\bF$-hyperbolic spaces with complex structures (when $\bF\neq\R$) cannot be totally geodesic due to a well-known result of H.Busemann \cite{Bu}, Theorem 47.4, p.331
(see also \cite{EI, AK}).
Nevertheless, as it was pointed out by G.D.Mostow in the complex hyperbolic geometry, such bisectors in $\bF$-hyperbolic spaces with complex structures are as close to being totally geodesic as possible. They are
minimal hypersurfaces of cohomogenity ($\dim_{\R} \bF - 1$), all
equivalent under $\bF$-hyperbolic isometries, and have a natural decomposition
into totally geodesic $\bF$-hypersurfaces ($\cong\fh {n-1}$), cf. \cite{A6}.

\subsection{Spherical CR-structures}

The $\F$-hyperbolic spaces $\fh m$ with complex structures (with $\bF\neq\R$) are naturally connected with spherical Cauchy-Riemannian (CR) structures at their infinity. Here a CR-structure on a real hypersurface in a complex manifold is defined by its largest subbundle in the tangent bundle that is invariant under the complex structure. In our case of the sphere $S^n=\p\fh m$, its spherical CR-structure is arising from describing this sphere as the boundary of the unit $\bF$-ball $B^m_{\bF}\subset \bF^m$. A CR-structure locally equivalent to the CR-structure of the boundary sphere of  the $\bF$-ball in $\bF^m$  is called a spherical CR-structure. A spherical CR-structure can also be viewed as a $(G,X)$-structure (see \cite{A3, A6}) with $X=S^n$ and $G$ as the isometry groups of one of $\bF$-hyperbolic spaces acting in the unit $\bF$-ball
$B^m_{\bF}\subset \bF^m$. A spherical CR-structure on a manifold/orbifold $M$ is said to be uniformizable if it is obtained from a discrete subgroup $\Ga\subset G$ by taking the quotient of the discontinuous action of $\Ga$ on an invariant connected component $\Om_0\subset \Om(\Ga)$ of the discontinuity set $\Om(\Ga)\subset \p\bF^m=S^n$, see \cite{A3, A6}.

\section{Discrete groups with $\La(G)=\p\fh n$}

Here we address a necessary condition for rigidity of deformations of locally symmetric $\R$-rank one manifolds/orbifolds $M$ universally covered by one of $\bF$-hyperbolic spaces $\fh n$ over division algebras $\bF$ - either
the real numbers $\br$, or
the complex numbers $\bc$, or
the quaternions $\bh$, or the Cayley numbers (octonions) $\bo$ (in the latter case $n=2$) -
the uniqueness of a spherical CR-structure (or
conformal one for real $\bF=\R$) at infinity $\fh n$ invariant for $\bF$-hyperbolic isometric action of the fundamental group
$\pi_1(M)\cong\Ga\subset\isom\fh n$ whose limit set $\La(\Ga)$ is the whole sphere at infinity $\p\fh n$.
As Theorem \ref{non-rigid} claims, such uniqueness of conformal/spherical CR-structures fails, and there are non-rigid discrete $\bF$-hyperbolic isometry groups $\Ga$ with $\La(\Ga)=\p\fh n$ and corresponding non-rigid $\bF$-hyperbolic manifolds/orbifolds $M=\fh n/\Ga$ (of infinite volume) without boundary at infinity.

\begin{pf}
In order to prove Theorem \ref{non-rigid}, here we provide a construction of desired non-rigid discrete groups
$\Ga\subset\isom\fh n$ with $\La(\Ga)=\p\fh n$, together with their non-trivial
deformations, i.e. non-trivial curves in the discrete faithful representation varieties of $\Ga$ passing through the inclusion $\Ga\subset\isom\fh n$. The discrete faithful representations $\rho_t\col\,\,\Ga\ra\isom\fh n$ defining such curves are obtained by families $\{\Phi_t\}$, $t\in(-\e,\e)$, of $\Ga$-equivariant
homeomorphisms with bounded distortion, $\Phi_t\col\ov{\fh n}\ra\ov{\fh n}$.

 We shall start with partitions of unit cubes either in the $(n-1)$-dimensional Euclidean space,
\begin{eqnarray}\label{cube}
Q=\{x\in\R^{n-1}\col|x_i|\leq 1/2\}\,,
\end{eqnarray}
or in the Carnot nilpotent group $\gtH=\bF^{n-1}\times\im\bF$ in the sphere at infinity $\p\fh n=\gtH\cup\{\infty\}$ whose 1-edges are parallel to the standard coordinate axes $\R\cdot e_k$, $k\geq 1$, in these spaces, and any two opposite faces are identified by corresponding mutually orthogonal unit left translations (in $\R^{n-1}$ or $\gtH$).
In other words, instead of the unit cube (\ref{cube}) in $\R^{n-1}$, in the nilpotent geometry in $\gtH$ we consider Carnot unit cube centered at the origin
$$
Q\subset\gtH\cong\F^{m-1}\times\im\F=\p\fh m\bs\{\infty\}
$$
generated by mutually orthogonal unit left translations in the nilpotent lattice (in $\gtH\cong\F^{m-1}\times\im\F$ where commutator of two mutually orthogonal horizontal unit translations (including the translation by $e_1=(1,0)\in\F^{m-1}\times\im\F$) is a vertical unit translation in the center $\{0\}\times\im\bF$ of the Carnot group $\gtH$.

The classical partition of the Euclidean 3-cube is known as the Sierpi\'{n}ski  carpet generalizing the Cantor set' partition of the unit segment. Our nilpotent Sierpi\'{n}ski carpet in $Q\subset\gtH$ defined by partitions of cubes in Carnot groups $\gtH$ should obey their nilpotent geometry.

For our (Carnot or Euclidean) unit cube $Q$ we consider an increasing sequence of odd integers $\{k_j\}$, $j\geq 1$.  Namely, in the first step, we split up the cube $Q$ into $k^{m-1}_1$ congruent subcubes $Q_{1i}$, where $m-1$ is the real dimension of the space, i.e. $m=n$ in the Eulidean case, and in the Carnot group $\gtH=\bF^{n-1}\times\im\bF$, $m=\dim_{\R}\fh n$.
The subcubes $Q_{1i}\subset Q$ create a tesselation of the cube $Q$. These subcubes are $\gtH$-isometric to each other under left translations given by the combinations (with integer coefficients) of our $(m-1)$ mutually orthogonal left translation vectors whose length is decreased by using nilpotent Carnot (or Euclidean) dilation $D_{\la}$ with $\la=1/k_1$, cf. (\ref{dilation}).

Then we cast out the interior of the central subcube $Q_{1i_0}$.  Similarly, in
the second step, we split up each of the remaining cubes $Q_{1i}\subset Q$, $i\neq i_0$, into
$k^{m-1}_2$ $\gtH$-isometric subcubes $Q_{2j}$ (creating tesselations of $Q_{1i}$)  and cast out the interiors of the central ones. We continue this process by using our increasing sequence of odd integers $\{k_m\}$.
Suitably choosing the increasing integers $k_j$, e.g. $k_j=3^j$, and
continuing the process of casting out, we finally obtain a limit continuum ($(m-1)$-dimensional Euclidean or Carnot Sierpi\'{n}ski carpet)
$K\subset Q$ whose Lebesgue $(m-1)$-measure $m_L(K)$ may equal any positive
number in the interval $[0,1)$, see \cite{H}. In the 3-dimensional Euclidean space, with the constant sequence $k_j=3$, the Lebesgue 2-measure $m_L(K)$ of the classical Sierpi\'{n}ski carpet $K$ is zero, and the Hausdorff dimension of this carpet is $\log 8/\log 3\approx 1.8928$.

In the closure of the half-space model of $\fh n$, $\ov{\fh n}\bs\{\infty\}$,
we consider the orthogonal projection $p_1$ to
the first coordinate axis $\R e_1$ in $\p\fh n$ (either the Euclidean space or the Carnot group $\bF^{n-1}\times\im\bF$) where in the latter case $\R e_1\subset\bF^{n-1}\times\im\bF$ is the first horizontal axis. Now let
$\Da_1=p_1(Q\bs K)$ be the subset in $[-1/2, 1/2]\subset\R e_1$ obtained as
the $p_1$-projection of the set $\Da\subset Q$, $\Da=Q\bs K$, where $\Da=Q\bs K$ consists of all casted out subcubes $Q_{lj_i}\subset Q$.
Then we define a map $f_t\col Q\times [0, \infty)\ra\ov{\fh n}\bs\{\infty\}$,
 $t>0$, as
\begin{eqnarray}\label{homeo}
f_t(x_1,\ldots,x_{k(n-1)},\ldots,x_{kn})=(\psi_t(x_1),x_2,\ldots,x_{k(n-1)},\ldots,x_{kn})\,,
\quad\psi_t(x_1)=\int^{x_1}_0
\vp_t(y)dy\,,
\end{eqnarray}
where $k=\dim_{\R}\bF$ and the measurable function $\vp_t(y)$ is
$$
\vp_t(y)=\begin{cases}
1\,, & \text{if}\quad y\in\Da_1\,;\\
t\,, & \text{if}\quad y\in\left[-\frac 12,\frac 12\right]\bs\Da_1\,.\end{cases}
$$

Such defined homeomorphism $f_t$ preserves (horizontal) subspaces $\bF^{n-1}\times\{v\}$ where its linear dilatation equal $t$ if $t\geq 1$, and it is $1/t$ if $0<t<1$.  The map $f_t$ transforms the qubic column
$Q\times [0, \infty)\subset\ov{\fh n}\bs\{\infty\}$ to the column
$f_t(Q\times [0, \infty))=Q_{t_1}\times [0, \infty)$, where the basis
box $Q_{t_1}\subset\bF^{n-1}\times\im\bF$ is
obtained from the cube $Q$ by its stretching/compressing along the first
coordinate axis $\R e_1$ by a factor $t_1=t_1(t)>0$.  This number $t_1$ depends only on $t$ and
the sequence $\{k_j\}$, and can be equal to $1$ only for $t=1$. We notice that the "vertical" edges of this
box $Q_{t_1}$ correspond to commutators of translations in $\gtH=\bF^{n-1}\times\im\bF$ along its horizontal edges.
For non-rigid initial lattice in $\gtH$ (with fundamental unit cube $Q$), these horizontal translations generate the (deformed) lattice in $\gtH$ whose fundamental polyhedron is the stretched box $Q_{t_1}$.
Also we notice that for each connected component of the  casted out set $\Da\subset Q$ (i.e. for a casted out subcube $Q_{li_j}\subset\Da\subset Q$), the restriction of the map
$f_t$ to $Q_{li_j}\times [0, \infty)$ is an $\bF$-hyperbolic isometry with unique fixed point at $\infty$ (a parabolic element in $\isom\fh n$) inducing a translation by a vector whose direction and length depend on the component $Q_{li_j}\subset\Da$, see Figure below.

Now we define isomorphic discrete $\F$-hyperbolic groups $\Ga,\Ga_t\subset\isom\fh m$. We need to create a dense covering $\Sa$ of the set $\Da\subset Q$ (all casted out subcubes $Q_{lj_i}\subset Q$) by disjoint closed balls $B_i\subset\Da\subset Q$, $i\in I$. In the real hyperbolic case, these balls $B_i$ are round balls. In the case of Hermitian hyperbolic structures $\bF\neq\R$, the balls $B_i\in\Sa$ are
disjoint Kor\'{a}nyi-Cygan closed balls $B_i\subset\Da\subset Q$ whose boundary spheres are (spinal) spheres at infinity of $\bF$-hyperbolic bisectors $\mathfrak S_i\subset\fh n$, see (\ref{bisector}), i.e. spheres $S_i$ with respect to the Kor\'{a}nyi-Cygan metric
$\rho_c$ in $\gtH=\bF^{n-1}\times\im\bF$, cf. (\ref{sphere}). The homeomorphism $f_t$ maps the family $\Sa=\{B_i\}$
to a new family $\Sa_t=\{B^t_i\}$ of disjoint closed balls in $\p\fh n\bs\{\infty\}$
which densely cover the image $f_t(\Da)\subset Q_{t_1}$.

\includegraphics{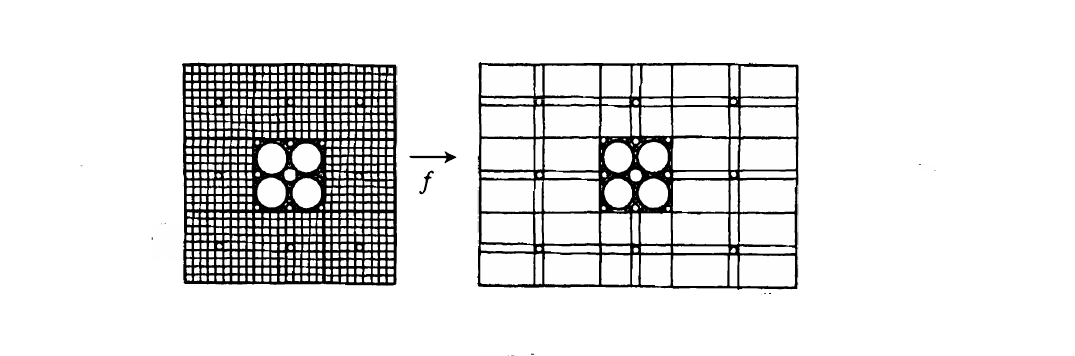}

We define discrete $\bF$-hyperbolic subgroups
$H$ and $H_t$ in $\isom\fh n$ (the half-space model of $\fh n$) as the groups generated by involutions in
boundary spheres $\{S_i=\p B_i\}$ and
$\{S^t_i=\p B^t_i\}$, respectively, i.e. the free products of order two subgroups generated by involutions $\sci_i$ in bisectors $\gtS_i\subset\fh n$, $\p\gtS_i=S_i$, see (\ref{bisector}). Also we define discrete parabolic subgroups $E$ and
$E_t\subset\isom\fh n$ acting in $\gtH=\p\fh n\bs\{\infty\}$ as lattices. In the real hyperbolic case, these groups are free Abelian rank $(n-1)$ groups whose generating Euclidean translations (acting in $\R^{n-1}$) identify the opposite faces of the cube $Q$ and the rectangular box $Q_{t_1}$, correspondingly. For Hermitian structures $\bF\neq\R$, these groups are lattices in the Carnot nilpotent
group $\gtH=\bF^{n-1}\times\im\bF$ generated by horizontal left translations in $\gtH$ identifying opposite faces of the cube $Q$ and the stretched box $Q_{t_1}$, correspondingly.

Now we can use the Klein Combination Theorem  \cite{AX, A3, A4, A6, FK}
to form the desired $\bF$-hyperbolic discrete groups $\Ga$ and $\Ga_t$ as the following free
products,
\begin{eqnarray}\label{groups}
\Ga=E\ast H\quad\text{and}\quad \Ga_t=E_t\ast H_t\,,
\end{eqnarray}
which act by isometries in the upper half-space model of the $\bF$-hyperbolic space $\fh n$.
Obviously, their limit sets $\La(\Ga)$ and $\La(\Ga_t)$ coincide with the whole sphere at infinity
$\p\fh n$. Denoting by $B_i^+$ the component of $\fh n\bs\gtS_i$ (the complement of a bisector $\gtS_i$, $S_i=\p\gtS_i$) having at its infinity the ball $B_i$, and by $B_i^+(t)$ the component of $\fh n\bs\gtS^t_i$ (the complement of a bisector $\gtS^t_i$, $S^t_i=\p\gtS^t_i$) having at its infinity the ball $B^t_i$,
we can choose fundamental $\bF$-hyperbolic
polyhedra $P(\Ga),P(\Ga_t)\subset\fh n$ for discrete groups $\Ga, \Ga_t\subset \isom\fh n$ as
\begin{eqnarray}\label{polyhedra}
P(\Ga)=(Q\times\R_+)\bs\bigcup_i B^+_i\,,\quad
P(\Ga_t)=(Q_{t_1}\times\R_+)\bs\bigcup_i B^+_i(t)\,.
\end{eqnarray}

Since the constructed in (\ref{homeo}) homeomorphisms $f_t$ provide mappings between polyhedra (\ref{polyhedra}), $f_t(P(\Ga))=P(\Ga_t)$, which are $\Ga$-equivariant on the boundary of
the polyhedron $P(\Ga)$ in $\fh n$, one can equivariantly extend the mappings $f_t\big\vert_{P(\Ga)}$ to
 $\Ga$-equivariant homeomorphisms $\Phi_t\col\,\fh n\ra\fh n$ whose
linear dilatations are the same as that of $f_t$. It clearly follows from the construction that
$\Ga_t=\Phi_t\Ga\Phi^{-1}_t=\Phi^*_t(\Ga)$.

Therefore the obtained continuous family $\Phi_t\col\,\fh n\ra\fh n$ of $\Ga$-equivariant homeomorphisms creates a continuous family of deformed discrete groups $\Ga_t=\Phi_t\Ga\Phi^{-1}_t$ isomorphic to the group $\Ga$ and defines a curve in the variety of conjugacy classes of discrete faitful representation of the group $\Ga$ passing through its inclusion $\Ga\subset\isom\fh n$.
This curve provides a non-trivial deformation of the inclusion $\Ga\subset\isom\fh n$ whose limit set is the whole sphere at infinity $\p\fh n$. Its non-triviality follows from an observation on the strict difference of the $\bF$-hyperbolic translation lengths of two isomorphic cyclic hyperbolic subgroups $A^t\subset\Phi^*_t(A)$ and $A^{t'}\subset\Phi^*_{t'}(A)$ with $t\neq t'$.  Here a cyclic hyperbolic subgroup $A\subset\Ga$ (with an invariant real hyperbolic geodesic in $\fh n$) is generated by involutions $\sci_i$ and $\sci_j$ in the bisectors $\gtS_i$ and $\gtS_j$ whose (spinal) spheres $S_i=\p B_i$ and $S_j=\p B_j$ at infinity are projected by the projection $p_1$ to disjoint connected components of the set $\Da_1=p_1(Q\bs K)$.
\end{pf}

\begin{rem}\label{contact-m} For Hermitian $\bF$-hyperbolic spaces $\fh n$ with $\bF\neq\R$,
we leave open the question whether the homeomorphism of the sphere at infinity $\p\fh n$ induced by the constructed $\Ga$-equivariant homeomorphism $\Phi_t\col\fh n\ra\fh n$ with bounded horizontal linear dilatation is a contactomorphism of the Carnot group at infinity
$\gtH=\p\fh n\bs\{\infty\}$. This is a necessary condition for its quasiconformality (which is unexpected in the quaternionic and octonionic hyperbolic spaces due to \cite{P}), cf. \cite{Mo2} and \cite{A5}.
\end{rem}

In the real hyperbolic space $\rh n$, the constructed $\Ga$-equivariant homeomorphisms
$\Phi_t\col\,\rh n\ra\rh n$ with bounded linear dilatations are extended to quasisymmetric $\Ga$-equivariant homeomorphisms of
$\ov{\rh n}=\rh n\cup\p\rh n$ whose restriction to the sphere at infinity $\p\rh n$ are unique quasiconformal self-homeomorphisms of $\p\rh n$ induced by the isomorphisms $\Phi^*_t\col\,\Ga\ra\Ga_t$ of discrete groups, see \cite{A3}.

Another approach to rigidity of locally symmetric rank one spaces $M$, or equivalently, of $\bF$-hyperbolic discrete holonomy groups $\pi^{\orb}_1(M)\cong\Ga\subset\isom\fh n$ whose limit sets coincide with the whole sphere at infinity $\p\fh n$, is based on Dennis Sullivan \cite{S} study of rigidity of real hyperbolic 3-manifolds. In this way, looking at dynamics of corresponding discrete holonomy groups $\Ga$ in the $\bF$-hyperbolic space $\fh n$ (see \cite{A3, A6}), one can connect rigidity of such groups and manifolds/orbifolds with zero
Lebesgue measure in the sphere at infinity, $m_{\infty}(\ov{P(\Ga)}\cap\p\fh n)=0$, of the limit subset on the boundary of the fundamental polyhedra $P(G)\subset\fh n$. Fixing a
point $0\in M$, we denote by $M(r)=\{x\in M\col\, d(0,x)\leq r\}$ the
$r$-neighborhood of this point in $M$ with respect to the $\bF$-hyperbolic
distance $d$, and $D(r)\subset\fh n$ a hyperbolic ball of radius $r>0$. Then one has:

\begin{rem}\label{dynamics}
The above condition could be characterized by any of the following equivalent
conditions, cf. \cite{A3, A6}:
\begin{enumerate}
\item{$\lim_{r\ra\infty}\vol M(r)/\vol D(r)=0$.}
\item{The action of the holonomy group $\Ga$ on the sphere $\p\fh n$
is conservative.}
\item{The horospherical limit set $\La_h(\Ga)$ has the full measure
in $\p\fh n$.}
\end{enumerate}
\end{rem}

\section{Deformations of $\bF$-hyperbolic manifolds trivial in the boundary at infinity}

A slight modification of our construction in the previous section addresses the existence of non-trivial deformations of $\bF$-hyperbolic manifolds trivial in their boundary at infinity. In other words, one is addressing deformations of the quotient spaces (Kleinian mani\-folds/orbifolds)
$M(\Ga)=[\fh n\cup\Om (\Ga)]/\Ga$ where $\Om(\Ga)\subset\p\fh n$ is the non-empty discontiunuity set of a discrete group
$\Ga\subset\isom\fh n$, with trivial restriction of such deformations to the boundary $\p M(\Ga)=\Om (\Ga)/\Ga$. This also answers another question
of L.~Bers \cite{B} and S.~Krushkal \cite{Kr1, Kr2} and its analogue for locally symmetric rank one manifolds with Hermitian structures, i.e. for $\bF$-hyperbolic structures with $\bF\neq\R$. In the latter case, the boundary of such Hermitian Kleinian manifolds/orbifolds $M(\Ga)$ (non-compact in general) has the natu\-ral spherical CR-structure covered by the spherical CR-structure on the discontinuity set $\Om(\Ga)\subset\p\fh n$. Moreover we should note that, in the case of a compact $M(\Ga)$ with Hermitian structure (i.e. for convex co-compact groups $\Ga\subset\isom\fh n$, see \cite{A3, A6}; in general, it is enough to have only one compact boundary component), an application of Kohn-Rossi \cite{KR} extension theorem shows that the boundary of $M(\Ga)$ (and the discontinuity set $\Om(\Ga))$ must be connected, and the
limit set $\La(G)$ is in some sense small. In these cases one has rigidity of deformations, cf. \cite{NR}, \cite{GM}, \cite{C}, \cite{GS}.
Our case is different - here we shall prove the non-rigidy result formulated in Theorem \ref{CR-boundary}.

\begin{pf} We can define desired non-rigid discrete $\bF$-hyperbolic groups $\Ga\subset\isom\fh n$ with non-empty discontinuity set
$\Om(\Ga)\subset\p\fh n$ as subgroups of the discrete groups constructed in the proof of Theorem \ref{non-rigid}, see (\ref{groups}).  Namely, in the definition of those discrete $\bF$-hyperbolic
groups (\ref{groups}), we adjust our dense ball-covering $\Sa$ of the set $\Da\subset Q$
 of removed subcubes $Q_{li_j}$ in the unit cube $Q$, see (\ref{cube}). It consists of disjoint closed balls $B_i\subset\Da\subset Q$, $i\in I$, whose boundary spheres are (spinal) spheres at infinity of $\bF$-hyperbolic bisectors $\gtS_i\subset\fh n$ (for $\bF\neq\R$, closed round balls with respect to the Kor\'{a}nyi-Cygan metric $\rho_c$ in the Carnot group $\gtH=\bF^{n-1}\times\im\bF$, see (\ref{metric})).
Here we just remove from the family $\Sa$ a ball $B_{i_0}$ (or a finite subset of balls $B_{i_1},\ldots,B_{i_k}$ from $\Sa$). Remaining balls in the family $\Sa$ provide remaining generators for discrete free groups $H$ and $H_t$, and we define our desired discrete groups $\Ga, \Ga_t\subset\isom\fh n$ as in (\ref{groups}). Then the discontinuity set $\Om(\Ga)\subset\p\fh n$ is the $\Ga$-orbit of the open ball $\rint B_{i_0}$ (or the finitely many open balls
$\rint B_{i_j}$, $j=1,\ldots,k$, at the sphere at infinity $\p\fh n$.

Upon this adjustment, our construction in the proof of Theorem \ref{non-rigid} gives us a new deformation $\{\Phi^0_t\}_{t\geq 0}$, $\Phi^0_0=\id$, of the new group $\Ga$ where
$\Ga$-equivariant homeomorphisms $\Phi^0_t\col\,\ov{\fh n}\ra\ov{\fh n}$
conjugate $\Ga$ to the discrete groups
$\Ga_t=\Phi^0_t\Ga(\Phi^0_t)^{-1}$. In the $\bF$-hyperbolic space $\fh n$, the $\Ga$-equivariant homeomorphisms $\Phi^0_t$ have the same linear dilatations as that of the constructed in (\ref{homeo}) homeomorphisms $f_t$.
Their restrictions $\Phi^0_t\big\vert_{\Om(\Ga)}$ to connected components of the discontinuity set
$\Om(Ga)\subset\p\fh n$ are restrictions of some $\bF$-hyperbolic isometries. Depending on connected components of $\Om(Ga)$ (in casted out subcubes $Q_{li_j}\subset Q$), these $\bF$-hyperbolic isometries are Euclidean translations
in $\R^{n-1}=\p\rh n\bs\{\infty\}$ if $\bF=\R$, and they are left Carnot translations in the Carnot group
$\gtH=\bF^{n-1}\times\im\bF$ for $\bF\neq\R$. Furthermore, this
deformation of the group $\Ga$ defined by our family of homeomorphisms $\{\Phi^0_t\}$ is non-trivial, and linear dilatations of used homeomorphisms
differ from 1 only on the limit set $\La(\Ga)$ (of positive Lebesgue measure on the sphere at infinity $\p\fh m$).

In the case of real hyperbolic spaces of dimension $n=2$ or $n=3$, that is for Kleinian groups on the extended complex
plane $\ov\C$, it means that the complex characteristics $\mu_t(z)$ of
our $\Ga$-equivariant quasiconformal homeomorphisms $\Phi^0_t\big\vert_{\ov\C}$ (solutions of the Beltrami equations) are
supported on the corresponding subset of the limit set $\La(\Ga)$ of positive Lebesgue $(n-1)$-measure.
\end{pf}

\begin{rem}\label{measure}
The above constructions in the proofs of
Theorems \ref{non-rigid}  and \ref{CR-boundary} of non-rigid discrete $\bF$-hyperbolic isometry groups $\Ga\subset\isom\fh n$ over division algebras $\bF=\R\,,\C\,,\bh\,,\text{or}\,\, \bo$ are based on their crucial property that the
subset of limit points $\La(\Ga)\cap\ov{P(\Ga)}$ on the boundary of the fundamental polyhedron
$P(\Ga)$ in (\ref{polyhedra}) has a positive Lebesgue measure in the sphere at infinity $\p\fh n$. This property directly follows from the positive Lebesgue measure $m_L(K)>0$ of the constructed Sierpi\'{n}ski carpet $K\subset Q$ in the unit cube $Q\subset\p\fh n\bs\{\infty\}$.
\end{rem}

\begin{rem}\label{convex hull}
Similarly to Remark \ref{dynamics}, one can relate the above Remark \ref{measure} non-rigidity condition to a dynamical condition on rigidity of deformations of locally $\bF$-hyperbolic mani\-folds/orbifolds $M$ where the asymptotics of $\bF$-hyperbolic volume of $M$ is replaced by the asymptotics of $\bF$-hyperbolic volume of its convex hull $C(M)\subset M$, cf. \cite{A3, A6}.
\end{rem}
\vfil\vfil

\bigskip

   Department of Mathematics, University of Oklahoma, Norman, OK 73019-0351, USA

     e-mail: apanasov\char`\@ ou.edu

     \end{document}